# Online decentralized tracking for nonlinear time-varying optimal power flow of coupled transmission–distribution grids

Wentian Lu[a,b*], Kaijun Xie[b], Mingbo Liu[b], Xiaogang Wang[a], Lefeng Cheng[a]

[a]School of Mechanical and Electric Engineering, Guangzhou University, Guangzhou 510006, China
[b]School of Electric Power Engineering, South China University of Technology, Guangzhou 510640, China

**Abstract**: The coordinated alternating current optimal power flow (ACOPF) for coupled transmission–distribution grids has become crucial to handle problems related to high penetration of renewable energy sources (RESs). However, obtaining all system details and solving ACOPF centrally is not feasible because of privacy concerns. Intermittent RESs and uncontrollable loads can swiftly change the operating condition of the power grid. Existing decentralized optimization methods can seldom track the optimal solutions of time-varying ACOPFs. Here, we propose an online decentralized optimization method to track the time-varying ACOPF of coupled transmission–distribution grids. First, the time-varying ACOPF problem is converted to a dynamic system based on Karush–Kuhn–Tucker conditions from the control perspective. Second, a prediction term denoted by the partial derivative with respect to time is developed to improve the tracking accuracy of the dynamic system. Third, a decentralized implementation for solving the dynamic system is designed based on only a few information exchanges with respect to boundary variables. Moreover, the proposed algorithm can be used to directly address nonlinear power flow equations without relying on convex relaxations or linearization techniques. Numerical test results reveal the effectiveness and fast-tracking performance of the proposed algorithm.

| **Nomenclature** | |
|---|---|
| *Sets and Indices* | |
| $N$ | Set of buses |
| $i$ | Index of bus |
| $N(i)$ | Set of buses connected to bus $i$ |
| $t$ | Index of time slot |
| $N_\text{g}$ | Index of conventional generator |
| $N_\text{res}$ | Index of RES |
| $\Gamma$ | Set of lines |
| $\Omega_\text{D}$ | Set of distribution systems |
| *Parameters* | |
| $P_i^\text{d}(t), Q_i^\text{d}(t)$ | Active and reactive power demand at bus $i$ at time $t$ |
| $P_i^\text{av}(t)$ | Maximum available active power output of the RES connected with bus $i$ |
| $\theta_i$ | Maximum allowable power-factor angle of the RES at bus $i$ |
| $S_i^\text{av}$ | Rated apparent power of the RES at bus $i$ |
| $\underline{V}_i, \overline{V}_i$ | Lower and upper bounds of voltage magnitude at bus $i$ |
| $\underline{P}_i^\text{g}, \overline{P}_i^\text{g}$ | Lower and upper active power output bounds of the generator at bus $i$ |
| $c_{\text{g}2,i}, c_{\text{g}1,i}, c_{\text{g}0,i}$ | Quadratic, linear, and constant terms in the quadratic generation cost function of the generator at bus $i$ |

---

[*] Corresponding author

E-mail: hnlgtiantian@163.com (Wentian Lu); ep_xiekj@mail.scut.edu.cn (Kaijun Xie); epmbliu@scut.edu.cn (Mingbo Liu); wxg@gzhu.edu.cn(Xiaogang Wang); chenglefeng@gzhu.edu.cn (Lefeng Cheng).

| $c_i^p, c_i^q$ | Quadratic and linear terms in the quadratic penalize function for the active power curtailment and reactive power generation of the RES at bus $i$ |
|---|---|
| $\bar{S}_{ij}$ | Maximum transmission capacity of the branch $ij$ |
| $g_{ij}, b_{ij}$ | Conductance and susceptance of the branch $ij$ |
| *Variables* | |
| $P_i^g(t), Q_i^g(t)$ | Real and reactive power outputs of the generator at bus $i$ at time $t$ |
| $P_i^{res}(t), Q_i^{res}(t)$ | Real and reactive power outputs of the RES at bus $i$ at time $t$ |
| $e_i(t), f_i(t)$ | Real and imaginary parts of the voltage at bus $i$ at time $t$ |
| $P_{ij}(t), Q_{ij}(t)$ | Active and reactive power of the branch $ij$ at time $t$ |

## 1. Introduction

### 1.1. Background

In current power grid operations, the transmission system (TS) and distribution systems (DSs) are physically connected but separately managed by the TS and DS operators, respectively. Each operator considers the other system as an equivalent; that is, the TS operator regards the DS as an equivalent load, and the DS operator regards the TS as an infinite equivalent generator [1]. With the increasing incorporation of renewable energy sources (RESs), such as photovoltaics (PVs) and wind turbines (WTs), into distribution systems, regarding DSs as equivalent loads in TS alternating current optimal power flow (ACOPF) computation causes numerous problems, such as power mismatch at boundary buses, voltage rise, and resource wastage [2]. Thus, increased coordination is required between the TS and DSs.

To solve this problem, a solution is required to centrally coordinate the TS and DSs. However, this approach requires full knowledge of the integrated TS and DSs, and solving such a large-scale ACOPF centrally is not feasible because of data privacy concerns. Therefore, a decentralized coordinated optimization method is preferable for solving coordinated transmission and distribution ACOPF. This solution can converge to a centralized solution with limited information exchange.

### 1.2. Related literature

From the coordination perspective, several studies have focused on decomposition approaches for decentralized operation problems. The algorithms can be categorized into three categories. The first set of methods is based on the Lagrangian relaxation (LR), in which a decomposable problem is constructed by relaxing the consistency constraint to the Lagrangian function. Subsequently, each operator solves its subproblem using Lagrange multipliers exchange among operators. These methods include classic LR [3], dynamic multiplier-based LR [4], and augmented Lagrangian decomposition (ALD) [5]. As an extension to the LR, in the ALD, a quadratic term about consistency constraint is added to enhance convergence. However, this addition increases decomposition difficulty. Extensive methods, including the auxiliary problem principle (APP) [6] with the application of multi-area optimal power flow (OPF) computation [7,8], the alternating direction method of multipliers (ADMM) [9] with its improvement on fully distributed ACOPF [10], and the analytical target cascading method for unit commitment problems [11,12], have been conducted on decomposing this quadratic term of ALD. However, the aforementioned algorithms typically require tuning of parameters in the multiplier updating process, which results in slow convergence or convergence failure.

The second set of decentralized optimization methods is based on the first-order Karush–Kuhn–Tucker (KKT) optimality condition decomposition (OCD). However, the method only features a first-order convergence rate. To enhance convergence, many methods, such as the approximate Newton direction (AND) method [13] with its multi-area OPF application [14], the Biskas algorithm for decentralized DC-OPF computation [15,16], and heterogeneous decomposition (HGD) for coordinated transmission and distribution optimal reactive power flow (ORPF) [17], economic dispatch [18], and OPF [19,20], have been proposed for decomposing the KKT condition. Although excellent computation accuracy and speed have been achieved, the convergence stability could be considerably affected by system size and system partitioning. Decentralized methods based on second-order KKT decomposition with a second-order convergence rate, such as the decomposition-coordination interior point method for multi-area ORPF [21] and distributed primal-dual interior point methods (PDIPM) for multi-area OPF [22], have been proposed. However, these methods require computation of the second-order Hessian inverse at each iteration, which could be prohibitively time-consuming for large-scale systems.

The last set of decentralized optimization algorithms for the transmission and distribution coordination are designed based on generating feasible cut constraints, such as generalized benders decomposition (GBD) [23,24], multiparametric programming [25–27], and the two-stage decomposition approach [28], which are used for decentralized models with a master–slave structure. However, linear cut constraints considerably influence the efficiency of these algorithms. Linearized or convexified techniques, such as second-order conic relaxations [23,26], semi-definite relaxations [24,29], a linearized AC power flow approach [25], a linear direct current OPF model [27], or a linearized DisFlow model [30–32], are required for the nonlinear ACOPF model. However, the optimal solution of the approximate model is not necessarily a feasible solution to the original ACOPF model.

With RESs being increasing incorporated in power grids, considering the uncertainty of RESs is critical in decentralized algorithms. Robust optimization (RO) is typically used to handle uncertainty in decentralized coordination of transmission and distribution grids, e.g., robust OPF [27,33]. Distributional RO is an improved version of RO [30,31]. In addition to RO, stochastic optimization [11,12], Monte Carlo simulations [34,35], and chance-constrained programming [36,37] are widely used uncertainty methods in power system decentralized optimization. However, these uncertainty methods generally bring an intractable computation burden and thwart online applications. A real-time optimal operation paradigm has emerged in decentralized algorithms; unlike uncertain models, short-term RES and load forecasts are accurate, and thus their prediction errors are assumed to be negligible in the real-time operating environment. In [38–41], online feedback-based distributed optimization was investigated in which local voltage measurements are implemented in real time and considerable information is shared with neighbors. However, a linearized model was adopted under the assumption that all the branches exhibited the same X/R ratio, which is not feasible in transmission networks. In [42], the update steps of the ADMM algorithm were modified to accommodate the voltage and branch flow measurements to achieve real-time distributed OPF. However, each RES communicates with all the other nodes and collect their measurements when updating its own power settings. This method is unsuitable for transmission and distribution coordination in terms of the communication structure.

*1.3. Contributions*

This study addressed the limitations of the aforementioned algorithms by proposing an online distributed optimization algorithm for timely tracking of the time-varying ACOPF for coupled transmission–distribution grids. Based on KKT conditions, the time-varying ACOPF problem is converted to a dynamic system from the control perspective. A prediction term, that is, the partial derivative with respect to time, is developed to improve the tracking accuracy for the dynamic system. A decentralized implementation structure is proposed based on only a few information exchanges with respect to boundary variables for solving the dynamic system. The major contributions of this study are summarized as follows:

1) An online decentralized tracking strategy is designed for coordinated ACOPF problems in coupled transmission–distribution grids. During each short sampling period, each DS only transforms a quadratic function with respect to the boundary variables' increment. Subsequently, the TS computes the boundary variables' increment and sends it to the corresponding distribution system. Therefore, the proposed decentralized implementation can protect the information privacy of each operator.

2) Unlike conventional iterative-convergence algorithms, the proposed online tracking algorithm only executes one iteration during each short sampling period, including one-time algebraic computing, algebraic updating, information exchanging, and command sending. These operations can save considerable time.

3) A PDIPM-based dynamic system is constructed from the control perspective to directly handle nonlinear AC power flow constraints by introducing time-varying slack variables and barrier parameters. Thus, using any approximate model based on convex relaxations or linearization techniques is avoided.

This article is organized as follows. In Section 2, a time-varying centralized TDOPF problem is formulated. In Section 3, first, a centralized tracking approach based on the PDIPM is proposed to solve the TDOPF problem. Next, the online distributed tracking algorithm is described in detail. In Section 4, simulation results on illustrative systems are presented. The conclusion is presented in Section 5.

## 2. TDOPF problem formulation

### 2.1. Primary TDOPF model

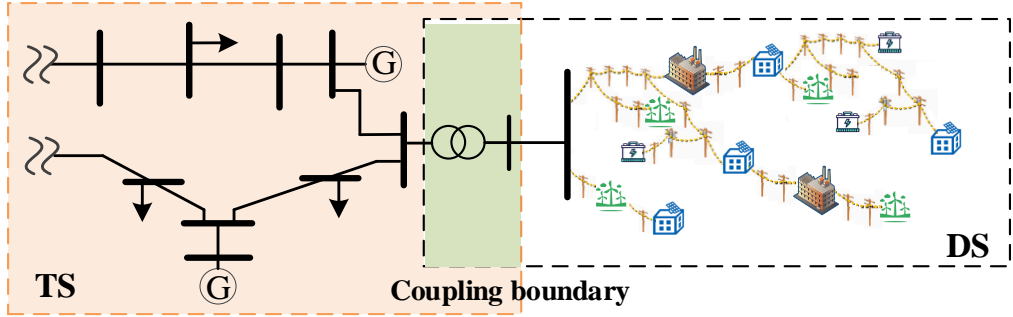

**Fig. 1. Topology of the integrated transmission and distribution system.**

A typical topology of the integrated transmission and distribution system is displayed in Fig. 1, in which TS and DS are connected by the tie line. This connection is defined as a coupling boundary. The TS and DS are physically connected but separately managed by the TS operator and DS operator, respectively. Assuming that the TS operator can access DS information, the optimal results for the integrated transmission and distribution system can be obtained by the following centralized transmission and distribution OPF (TDOPF) model:

$$\min \sum_{i \in N_g} C_i^g(P_i^g) + \sum_{i \in N_{res}} C_i^{res}(P_i^{res}, Q_i^{res}), \tag{1a}$$

subjective to
$$P_i^g + P_i^{res} - P_i^d(t) - \sum_{j \in N(i)} P_{ij} = 0, \forall i \in N, \tag{1b}$$

$$Q_i^g + Q_i^{res} - Q_i^d(t) - \sum_{j \in N(i)} Q_{ij} = 0, \forall i \in N, \tag{1c}$$

$$P_{ij} = g_{ij}(e_i^2 + f_i^2) - g_{ij}(e_i e_j + f_i f_j) + b_{ij}(e_i f_j - e_j f_i), ij \in \Gamma, \tag{1d}$$

$$Q_{ij} = -b_{ij}(e_i^2 + f_i^2) + b_{ij}(e_i e_j + f_i f_j) + g_{ij}(e_i f_j - e_j f_i), ij \in \Gamma, \tag{1e}$$

$$\underline{P}_i^g \leq P_i^g \leq \overline{P}_i^g, \forall i \in N_g, \tag{1f}$$

$$0 \leq P_i^{res} \leq P_i^{av}(t), \forall i \in N_{res}, \tag{1g}$$

$$\left|Q_i^{res}\right| \leq P_i^{res} \tan \theta_i, \forall i \in N_{res}, \tag{1h}$$

$$(P_i^{res})^2 + (Q_i^{res})^2 \leq (S_i^{av})^2, \forall i \in N_{res}, \tag{1i}$$

$$\underline{V}_i^2 \leq e_i^2 + f_i^2 \leq \overline{V}_i^2, \forall i \in N, \tag{1j}$$

$$P_{ij}^2 + Q_{ij}^2 \leq (\overline{S}_{ij})^2, \forall ij \in \Gamma, \tag{1k}$$

where $t$ is a continuous index of time. In (1a), the first term represents the generation costs of all conventional generators (e.g., diesels or turbines) and is expressed as the quadratic function active power: $C_i^g(P_i^g) = c_{g2,i}(P_i^g)^2 + c_{g1,i}P_i^g + c_{g0,i}$. For the second term, a quadratic function $C_i^{res}(P_i^{res}, Q_i^{res}) = c_i^p(P_i^{av} - P_i^{res})^2 + c_i^q(Q_i^{res})^2$ is used to minimize both the amount of active power curtailment and reactive power injected or absorbed of RES. Constraints (1b)–(1e) denote nonlinear AC power flow equations. Constraints (1f)–(1i) are the operational constraints of conventional generators and RESs. Constraints (1j) and (1k) describe the voltage magnitude limits of the buses and the transmission capacity limits of the lines, respectively.

Constraints (1b), (1c), and (1g) contain three time-varying parameters: the active and reactive power demand $P_i^d(t)$ and $Q_i^d(t)$, and the maximum available active power output of RESs $P_i^{av}(t)$. Thus, the aforementioned TDOPF model (1a)–(1k) is essentially a time-varying optimization problem.

*2.2. Reformulation of the TDOPF model*

The centralized TDOPF model (1a)–(1k) can be represented in a compact form to introduce the decentralized solution.

$$\min C_T(x_T^I(t), x^B(t)) + \sum_{k \in \Omega_D} C_{D,k}(x_{D,k}^I(t), x^B(t)), \tag{2a}$$

subjective to
$$g_T(x_T^I(t), x^B(t)) = 0, \tag{2b}$$
$$\bar{h}_T \leq h_T(x_T^I(t), x^B(t)) \leq \underline{h}_T, \tag{2c}$$
$$g_{D,k}(x_{D,k}^I(t), x^B(t)) = 0, k \in \Omega_D, \tag{2d}$$
$$\bar{h}_{D,k} \leq h_{D,k}(x_{D,k}^I(t), x^B(t)) \leq \underline{h}_{D,k}, k \in \Omega_D, \tag{2e}$$

where $x_T^I(t)$ represents the vector of the independent variable of the TS model, that is, variables that are only in the TS objective and constraints. Similarly, $x_{D,k}^I(t)$ is the vector of the independent variable of the $k^{th}$ DS model; TS and DS models are linked through the vector of boundary variables $x^B(t)$, which is specified by $x^B(t) = \{x_{T,k}^B(t) \bigcup x_{D,k}^B(t), k \in \Omega_D\}$. The detailed variable partition for the TDOPF model is displayed in Fig. 2.

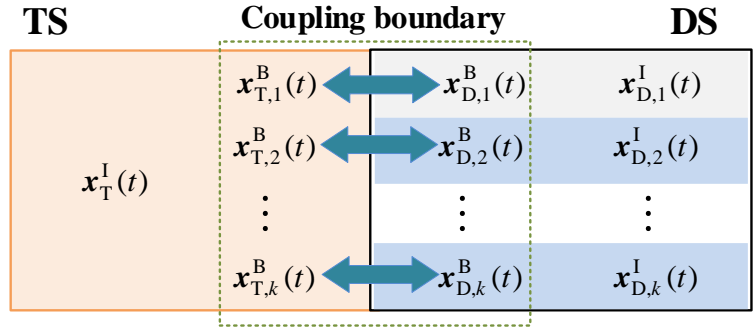

**Fig. 2. Variable partition for the TDOPF model.**

For exposition, we decomposed the centralized model (2a)–(2e) into a decentralized form, that is, a TS-OPF and DS-OPFs. For the TS, we have the following expression:

$$\min C_T(x_T^I(t)), \tag{3a}$$
$$\text{s.t. } g_T(x_T^I(t), x^B(t)) = 0, \tag{3b}$$
$$\bar{h}_T \leq h_T(x_T^I(t)) \leq \underline{h}_T, \tag{3c}$$

and for the DS $k$, $k \in \Omega_D$, we have the following expression:

$$\min C_{D,k}(x_{D,k}^I(t)), \tag{4a}$$
$$\text{s.t. } g_{D,k}(x_{D,k}^I(t), x^B(t)) = 0, \tag{4b}$$
$$\bar{h}_{D,k} \leq h_{D,k}(x_{D,k}^I(t)) \leq \underline{h}_{D,k}. \tag{4c}$$

From (3a)–(3c) and (4a)–(4c), because boundary variables exist between the TS and DSs, models (3a)–(3c) and (4a)–(4c) cannot be solved independently even for a given time instant. Furthermore, the fast-varying parameters ($P_i^{av}(t), P_i^d(t)$ and $Q_i^d(t)$) increase the difficulty of solvers for the decentralized TDOPF model.

## 3. Proposed method

The polynomial-time PDIPM in [43] is widely used for nonlinear optimization problems, especially for dealing with the nonlinear AC power flow constraints in power grids. However, the method is a centralized and iterative-convergence solver and cannot be used in an online application (such as a very fast timescale: i.e., 20 ms). In this study, we recast the PDIPM to fast track a KKT trajectory of the nonlinear time-varying models (3a)–(3c) and (4a)–(4c) in a decentralized manner.

*3.1. Building a dynamic system based on the KKT trajectory*

To simplify the description, we first unify (3a)–(3c) and (4a)–(4c) as follows:

$$\min f(x, t), \tag{5a}$$
$$\text{s.t. } g(x, t) = 0, \tag{5b}$$
$$\underline{h} \leq h(x, t) \leq \bar{h}. \tag{5c}$$

Next, we can construct the following Lagrange function:

$$L = f(\mathbf{x},t) - \mathbf{y}^T\mathbf{g}(\mathbf{x},t) - \mathbf{w}^T[\mathbf{h}(\mathbf{x},t)+\mathbf{u}-\overline{\mathbf{h}}] - \mathbf{z}^T[\mathbf{h}(\mathbf{x},t)-\mathbf{l}-\underline{\mathbf{h}}] - \mu(\sum \ln u_i + \sum \ln l_i), \quad (6)$$

where $\mathbf{y}$, $\mathbf{w}$, and $\mathbf{z}$ are the Lagrange multipliers (dual variables) for equality and inequality constraints, respectively, which also belong to time-varying variables; $\mathbf{u}$ and $\mathbf{l}$ are time-varying slack variables; $\mu$ is the barrier parameter with $\mu > 0$. Thus, the stationary point of (6) is the optimal solution of models (5a)–(5c), which satisfies the first-order optimal conditions, that is, KKT conditions. For simplicity, the aggregate variable vector is defined as $\boldsymbol{\lambda} := [\mathbf{x};\mathbf{y};\mathbf{w};\mathbf{z};\mathbf{u};\mathbf{l}]$, and the optimal solution as $\boldsymbol{\lambda}^*(t)=[\mathbf{x}^*(t);\mathbf{y}^*(t);\mathbf{w}^*(t);\mathbf{z}^*(t);\mathbf{u}^*(t);\mathbf{l}^*(t)]$. Therefore, we list the KKT conditions of models (5a)–(5c) in the following condensed form:

$$\nabla_\lambda L(\boldsymbol{\lambda}(t),t) = 0, \forall t \geq 0. \quad (7)$$

Applying the Newton's method to solve (7), the following dynamic system can be obtained:

$$\dot{\boldsymbol{\lambda}}(t) = -\left[\nabla_{\lambda\lambda}L(\boldsymbol{\lambda}(t),t)\right]^{-1}\nabla_\lambda L(\boldsymbol{\lambda}(t),t), \forall t \geq 0. \quad (8)$$

The trajectory $\boldsymbol{\lambda}(t)$ generated by (8) can approach a neighborhood around $\boldsymbol{\lambda}^*(t)$; however, this method cannot converge exactly to $\boldsymbol{\lambda}^*(t)$ because the problem is changing over time and so is the solution. The optimal solution $\boldsymbol{\lambda}^*(t)$ satisfies (7) for all times $t \geq 0$. Thus, the partial derivative of $\nabla_\lambda L(\boldsymbol{\lambda}^*,t)$ with respect to $t \in [0,T]$ should also be equal to $\mathbf{0}$, that is,

$$\mathbf{0} = \nabla_{\lambda\lambda}L(\boldsymbol{\lambda}^*(t),t)\dot{\boldsymbol{\lambda}}^*(t) + \nabla_{\lambda t}L(\boldsymbol{\lambda}^*(t),t), \forall t \geq 0, \quad (9)$$

and solving (9) for $\boldsymbol{\lambda}^*(t)$ yields another dynamic system as follows:

$$\dot{\boldsymbol{\lambda}}^*(t) = -\left[\nabla_{\lambda\lambda}L(\boldsymbol{\lambda}^*(t),t)\right]^{-1}\nabla_{\lambda t}L(\boldsymbol{\lambda}^*(t),t), \forall t \geq 0. \quad (10)$$

If the system parameters (e.g., the active and reactive power demand, and the maximal available active power output of RESs) are time-invariable, that is, $\nabla_{\lambda t}L(\boldsymbol{\lambda}^*(t),t)=0$, the right-hand side of the system (10) is equal to zero; otherwise, the system (10) can be considered a prediction term for correcting the tracking error of the dynamic system (8). Thus, a novel dynamic system can be constructed as follows:

$$\dot{\boldsymbol{\lambda}}(t) = -\left[\nabla_{\lambda\lambda}L(\boldsymbol{\lambda}(t),t)\right]^{-1}\left(\alpha\nabla_\lambda L(\boldsymbol{\lambda}(t),t) + \nabla_{\lambda t}L(\boldsymbol{\lambda}(t),t)\right), \quad (11)$$

where $\alpha$ is a positive weight coefficient. In the dynamic system (11), the first term from (8) can push $\boldsymbol{\lambda}(t)$ toward the optimum, and the second term from (10) can predict how the optimal solution changes over time by considering time variations of problem parameters. The general principle of building a dynamic system is displayed in Fig. 3, by which solving the nonlinear time-varying model (5a)–(5c) is converted into tracking a KKT trajectory.

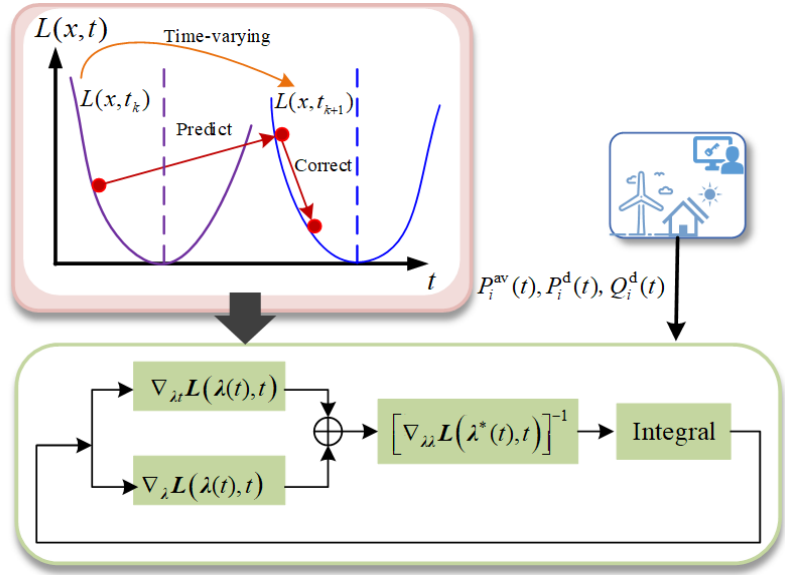

Fig. 3. Process of building a dynamic system for KKT conditions.

The general principle of developing a dynamic system can be extended to solving models (3a)–(3c) and (4a)–(4c). According to the superposition principle, a global variable vector can be defined as the union of the variable vector from TS and DSs, that is, $\hat{\lambda}(t) = \{\lambda_T(t) \cup \lambda_{D,k}(t), k \in \Omega_D\}$ (the superscript $\hat{}$ represents the union of variables from TS and DSs). The following centralized dynamic system can then be obtained:

$$\dot{\hat{\lambda}}(t) = -\left[\nabla_{\hat{\lambda}\hat{\lambda}} L(\hat{\lambda}(t),t)\right]^{-1} \left(\alpha \nabla_{\hat{\lambda}} L(\hat{\lambda}(t),t) + \nabla_{\hat{\lambda}t} L(\hat{\lambda}(t),t)\right), \quad (12)$$

where

$$\nabla_{\hat{\lambda}\hat{\lambda}} L(\hat{\lambda}(t),t) = \sum_{k \in \Omega_D} (\bar{E}_k)^T \left(\nabla_{\lambda_{D,k}\lambda_{D,k}} L_{D,k}(\lambda_{D,k}(t),t)\right) \bar{E}_k + (\bar{E}_T)^T \left(\nabla_{\lambda_T \lambda_T} L_T(\lambda_T(t),t)\right) \bar{E}_T, \quad (13a)$$

$$\nabla_{\hat{\lambda}} L(\hat{\lambda}(t),t) = \sum_{k \in \Omega_D} (\bar{E}_k)^T \left(\nabla_{\lambda_{D,k}} L_{D,k}(\lambda_{D,k}(t),t)\right) + (\bar{E}_T)^T \left(\nabla_{\lambda_T} L_T(\lambda_T(t),t)\right), \quad (13b)$$

$$\nabla_{\hat{\lambda}t} L(\hat{\lambda}(t),t) = \sum_{k \in \Omega_D} (\bar{E}_k)^T \left(\nabla_{\lambda_{D,k}t} L_{D,k}(\lambda_{D,k}(t),t)\right) + (\bar{E}_T)^T \left(\nabla_{\lambda_T t} L_T(\lambda_T(t),t)\right), \quad (13c)$$

where $\bar{E}_T$ is a $R_T \times R$ incidence matrix obtained from an $R \times R$ identity matrix by reserving rows with indices of $\lambda_T(t)$ in $\hat{\lambda}(t)$, that is, $\lambda_T(t) = \bar{E}_T \hat{\lambda}(t)$, and $R_T$ and $R$ denotes the number of elements in $\lambda_T(t)$ and $\hat{\lambda}(t)$, respectively. Similarly, $\bar{E}_k$ is a $R_k \times R$ incidence matrix obtained from an $R \times R$ identity matrix by reserving rows with indices of $\lambda_{D,k}(t)$ in $\hat{\lambda}(t)$, and $R_k$ denotes the number of elements in $\lambda_{D,k}(t)$, that is, $\lambda_{D,k}(t) = \bar{E}_k \hat{\lambda}(t), k \in \Omega_D$.

However, (12) is a centralized solution because its execution requires computing the Hessian inverse $\left[\nabla_{\hat{\lambda}\hat{\lambda}} L(\hat{\lambda}(t),t)\right]^{-1}$, and implementing such a global computation is not feasible for transmission and DSs because of privacy concerns. In the following subsection, we focus on solving the dynamic systems (13a)–(13c) in a decentralized fashion, that is, a protocol such that each DS only requires a few boundary information exchanges with the TS.

*3.2. Solving the dynamic system in a decentralized manner*

In theory, the KKT optimality trajectory can be tracked with a vanishing tracking error if the dynamic system (12) can be solved precisely. However, in practice, a suitable numerical method can be used to solve the dynamic system to obtain an approximate KKT trajectory. A computationally frugal choice is using the forward Euler method with a constant sampling period $\tau$ as follows:

$$\hat{\lambda}(t+\tau) = \hat{\lambda}(t) + \tau \Delta \hat{\lambda}(t), \forall t \geq 0, \quad (14)$$

where $\Delta \hat{\lambda}(t)$ is the variables' increments, that is, the right-hand side of the dynamic system (12). Thus, the variables' increments can be obtained by solving the following correction equation:

$$\nabla_{\hat{\lambda}\hat{\lambda}} L(\hat{\lambda}(t),t) \Delta \hat{\lambda}(t) = -\left(\alpha \nabla_{\hat{\lambda}} L(\hat{\lambda}(t),t) + \nabla_{\hat{\lambda}t} L(\hat{\lambda}(t),t)\right). \quad (15)$$

The decentralized implementation for solving (15) is based on only a few information exchanges with respect to boundary variables $x^B(t)$, which is displayed in Fig. 4.

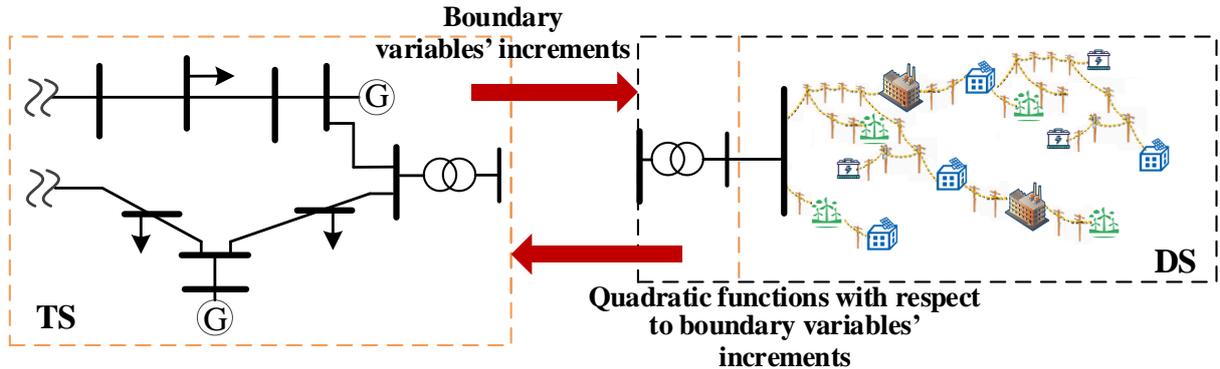

**Fig. 4. Framework based on only a few information exchanges with respect to boundary variables.**

*3.2.1 Quadratic function derivation*

Considering distribution system $k$ ($k \in \Omega_D$) as an example, on basis of the PDIPM, a sparse linear system (16) can be constructed with respect to the variable increments $\Delta x_{D,k}(t)$ and $\Delta y_{D,k}(t)$ as follows:

$$\begin{bmatrix} \boldsymbol{H}_{\mathrm{D},k} & \nabla_x^{\mathrm{T}} \boldsymbol{g}_{\mathrm{D},k}(\boldsymbol{x},t) \\ \nabla_x \boldsymbol{g}_{\mathrm{D},k}(\boldsymbol{x},t) & 0 \end{bmatrix} \begin{bmatrix} \Delta \boldsymbol{x}_{\mathrm{D},k}(t) \\ \Delta \boldsymbol{y}_{\mathrm{D},k}(t) \end{bmatrix} = \begin{bmatrix} \boldsymbol{R}_{k,x} \\ \boldsymbol{R}_{k,y} \end{bmatrix}, \ \forall k \in \Omega_{\mathrm{D}}. \tag{16}$$

After obtaining $\Delta \boldsymbol{x}_{\mathrm{D},k}(t)$ and $\Delta \boldsymbol{y}_{\mathrm{D},k}(t)$, the remaining $\Delta \boldsymbol{w}_{\mathrm{D},k}(t), \Delta \boldsymbol{z}_{\mathrm{D},k}(t), \Delta \boldsymbol{u}_{\mathrm{D},k}(t)$, and $\Delta \boldsymbol{l}_{\mathrm{D},k}(t)$ can be calculated as follows:

$$\Delta \boldsymbol{u}_{\mathrm{D},k}(t) = -(\boldsymbol{L}_{k,w} + \boldsymbol{L}_{k,wt}) - \nabla_x^{\mathrm{T}} \boldsymbol{h}_{\mathrm{D},k}(\boldsymbol{x}_{\mathrm{D},k},t) \Delta \boldsymbol{x}_{\mathrm{D},k}(t), \ \forall k \in \Omega_{\mathrm{D}}, \tag{17a}$$

$$\Delta \boldsymbol{l}_{\mathrm{D},k}(t) = (\boldsymbol{L}_{k,z} + \boldsymbol{L}_{k,zt}) + \nabla_x^{\mathrm{T}} \boldsymbol{h}_{\mathrm{D},k}(\boldsymbol{x}_{\mathrm{D},k},t) \Delta \boldsymbol{x}_{\mathrm{D},k}(t), \ \forall k \in \Omega_{\mathrm{D}}, \tag{17b}$$

$$\Delta \boldsymbol{w}_{\mathrm{D},k}(t) = -\boldsymbol{I}_{k,U}^{-1}(\boldsymbol{L}_{k,u} + \boldsymbol{L}_{k,ut}) - \boldsymbol{I}_{k,U}^{-1} \boldsymbol{I}_{k,W} \Delta \boldsymbol{u}_{\mathrm{D},k}(t), \ \forall k \in \Omega_{\mathrm{D}}, \tag{17c}$$

$$\Delta \boldsymbol{z}_{\mathrm{D},k}(t) = -\boldsymbol{I}_{k,L}^{-1}(\boldsymbol{L}_{k,l} + \boldsymbol{L}_{k,lt}) - \boldsymbol{I}_{k,L}^{-1} \boldsymbol{I}_{k,Z} \Delta \boldsymbol{l}_{\mathrm{D},k}(t), \ \forall k \in \Omega_{\mathrm{D}}, \tag{17d}$$

where

$$\boldsymbol{H}_{\mathrm{D},k} = -\left[ \nabla_{xx} f_{\mathrm{D},k}(\boldsymbol{x}_{\mathrm{D},k},t) - \nabla_{xx} \boldsymbol{g}_{\mathrm{D},k}(\boldsymbol{x}_{\mathrm{D},k},t) \boldsymbol{y}_{\mathrm{D},k}(t) - \nabla_{xx} \boldsymbol{h}_{\mathrm{D},k}(\boldsymbol{x}_{\mathrm{D},k},t)(\boldsymbol{w}_{\mathrm{D},k}(t) + \boldsymbol{z}_{\mathrm{D},k}(t)) \right]$$
$$- \nabla_x^{\mathrm{T}} \boldsymbol{h}_{\mathrm{D},k}(\boldsymbol{x}_{\mathrm{D},k},t)(\boldsymbol{I}_{k,L}^{-1} \boldsymbol{I}_{k,Z} - \boldsymbol{I}_{k,U}^{-1} \boldsymbol{I}_{k,W}) \nabla_x \boldsymbol{h}_{\mathrm{D},k}(\boldsymbol{x}_{\mathrm{D},k},t), \tag{18a}$$

$$\boldsymbol{R}_{k,x} = (\boldsymbol{L}_{k,x} + \boldsymbol{L}_{k,xt}) + \nabla_x^{\mathrm{T}} \boldsymbol{h}_{\mathrm{D},k}(\boldsymbol{x}_{\mathrm{D},k},t) \left[ \boldsymbol{I}_{k,L}^{-1}(\boldsymbol{L}_{k,l} + \boldsymbol{L}_{k,lt} + \boldsymbol{I}_{k,Z}(\boldsymbol{L}_{k,z} + \boldsymbol{L}_{k,zt})) + \boldsymbol{I}_{k,U}^{-1}(\boldsymbol{L}_{k,u} + \boldsymbol{L}_{k,ut} - \boldsymbol{I}_{k,W}(\boldsymbol{L}_{k,w} + \boldsymbol{L}_{k,wt})) \right], \tag{18b}$$

$$\boldsymbol{R}_{k,y} = -(\boldsymbol{L}_{k,y} + \boldsymbol{L}_{k,yt}). \tag{18c}$$

In (16)–(18c), $\boldsymbol{H}_{\mathrm{D},k}$ is the Hessian matrix; $\boldsymbol{I}_{k,U}, \boldsymbol{I}_{k,W}, \boldsymbol{I}_{k,L}$, and $\boldsymbol{I}_{k,Z}$ are diagonal matrices of $u$, $w$, $l$, and $z$; $\boldsymbol{L}_{k,x}, \boldsymbol{L}_{k,y}, \boldsymbol{L}_{k,w}, \boldsymbol{L}_{k,z}, \boldsymbol{L}_{k,u}$ and $\boldsymbol{L}_{k,l}$ denote the residuals of the KKT condition; and $\boldsymbol{L}_{k,xt}, \boldsymbol{L}_{k,yt}, \boldsymbol{L}_{k,wt}, \boldsymbol{L}_{k,zt}, \boldsymbol{L}_{k,ut}$ and $\boldsymbol{L}_{k,lt}$ denote the partial derivative of the residuals with respect to time $t$.

From (16), the second partial derivative of the residuals of the KKT condition with respect to time $t$ only exists on the left side of the equation. Thus, when the time-varying input parameters update, only the residual vectors $\boldsymbol{R}_{k,x}$ and $\boldsymbol{R}_{k,y}$ should be updated.

Solving the sparse linear equation (16) is equivalent to solving the following quadratic programming problem:

$$\min_{\Delta \boldsymbol{x}_{\mathrm{D},k}(t), \Delta \boldsymbol{y}_{\mathrm{D},k}(t)} \frac{1}{2} \Delta \boldsymbol{x}_{\mathrm{D},k}^{\mathrm{T}}(t) \boldsymbol{H}_{\mathrm{D},k} \Delta \boldsymbol{x}_{\mathrm{D},k}(t) - \boldsymbol{R}_{k,x}^{\mathrm{T}} \Delta \boldsymbol{x}_{\mathrm{D},k}(t) + \Delta \boldsymbol{y}_{\mathrm{D},k}(t)^{\mathrm{T}} (\nabla_x \boldsymbol{g}_{\mathrm{D},k}(\boldsymbol{x},t) \Delta \boldsymbol{x}_{\mathrm{D},k}(t) - \boldsymbol{R}_{k,y}). \tag{19}$$

According to the variable partition displayed in Fig. 2, rearranging $\Delta \boldsymbol{x}_{\mathrm{D},k}(t)$ and the corresponding matrices yields the following results:

$$\Delta \boldsymbol{x}_{\mathrm{D},k}(t) = \begin{bmatrix} \Delta \boldsymbol{x}_{\mathrm{D},k}^{\mathrm{I}}(t) \\ \Delta \boldsymbol{x}_{\mathrm{D},k}^{\mathrm{B}}(t) \end{bmatrix}, \tag{20a}$$

$$\boldsymbol{H}_{\mathrm{D},k} = \begin{bmatrix} \boldsymbol{H}_{\mathrm{D},k}^{\mathrm{II}} & \boldsymbol{H}_{\mathrm{D},k}^{\mathrm{IB}} \\ \boldsymbol{H}_{\mathrm{D},k}^{\mathrm{BI}} & \boldsymbol{H}_{\mathrm{D},k}^{\mathrm{BB}} \end{bmatrix}, \tag{20b}$$

$$\nabla_x \boldsymbol{g}_{\mathrm{D},k} = \begin{bmatrix} \nabla_x \boldsymbol{g}_{\mathrm{D},k}^{\mathrm{I}} & \nabla_x \boldsymbol{g}_{\mathrm{D},k}^{\mathrm{B}} \end{bmatrix}, \tag{20c}$$

$$\boldsymbol{R}_{k,x} = \begin{bmatrix} \boldsymbol{R}_{k,x}^{\mathrm{I}} & \boldsymbol{R}_{k,x}^{\mathrm{B}} \end{bmatrix}. \tag{20d}$$

Substituting (20a)–(20d) into (19), and considering $\Delta \boldsymbol{x}_{\mathrm{D},k}^{\mathrm{B}}(t)$ as a parameter, we can obtain the following parametric quadratic programming model:

$$\min_{\Delta \boldsymbol{x}_{\mathrm{D},k}^{\mathrm{I}}(t), \Delta \boldsymbol{y}_{\mathrm{D},k}(t)} \frac{1}{2} \left( \Delta \boldsymbol{x}_{\mathrm{D},k}^{\mathrm{I}}(t) \right)^{\mathrm{T}} \boldsymbol{H}_{\mathrm{D},k}^{\mathrm{II}} \Delta \boldsymbol{x}_{\mathrm{D},k}^{\mathrm{I}}(t) + \left( \Delta \boldsymbol{x}_{\mathrm{D},k}^{\mathrm{I}}(t) \right)^{\mathrm{T}} \left( \boldsymbol{H}_{\mathrm{D},k}^{\mathrm{IB}} \Delta \boldsymbol{x}_{\mathrm{D},k}^{\mathrm{B}}(t) - \boldsymbol{R}_{k,x}^{\mathrm{I}} \right)$$
$$+ \Delta \boldsymbol{y}_{\mathrm{D},k}(t)^{\mathrm{T}} \left( \nabla_x \boldsymbol{g}_{\mathrm{D},k}^{\mathrm{I}} \Delta \boldsymbol{x}_{\mathrm{D},k}^{\mathrm{I}}(t) - (\boldsymbol{R}_{k,y} - \nabla_x \boldsymbol{g}_{\mathrm{D},k}^{\mathrm{B}} \Delta \boldsymbol{x}_{\mathrm{D},k}^{\mathrm{B}}(t)) \right) + \frac{1}{2} \left( \Delta \boldsymbol{x}_{\mathrm{D},k}^{\mathrm{B}}(t) \right)^{\mathrm{T}} \boldsymbol{H}_{\mathrm{D},k}^{\mathrm{BB}} \Delta \boldsymbol{x}_{\mathrm{D},k}^{\mathrm{B}}(t) - \left( \Delta \boldsymbol{x}_{\mathrm{D},k}^{\mathrm{B}}(t) \right)^{\mathrm{T}} \boldsymbol{R}_{k,x}^{\mathrm{B}}. \tag{21}$$

Therefore, $\Delta \boldsymbol{x}_{\mathrm{D},k}^{\mathrm{I}}(t)$ and $\Delta \boldsymbol{y}_{\mathrm{D},k}(t)$ can be obtained by the following:

$$\begin{bmatrix} \Delta \boldsymbol{x}_{\mathrm{D},k}^{\mathrm{I}}(t) \\ \Delta \boldsymbol{y}_{\mathrm{D},k}(t) \end{bmatrix} = \begin{bmatrix} \boldsymbol{H}_{\mathrm{D},k}^{\mathrm{II}} & (\nabla_x \boldsymbol{g}_{\mathrm{D},k}^{\mathrm{I}})^{\mathrm{T}} \\ \nabla_x \boldsymbol{g}_{\mathrm{D},k}^{\mathrm{I}} & 0 \end{bmatrix}^{-1} \begin{bmatrix} \boldsymbol{R}_{k,x}^{\mathrm{I}} - \boldsymbol{H}_{\mathrm{D},k}^{\mathrm{IB}} \Delta \boldsymbol{x}_{\mathrm{D},k}^{\mathrm{B}}(t) \\ \boldsymbol{R}_{k,y} - \nabla_x \boldsymbol{g}_{\mathrm{D},k}^{\mathrm{B}} \Delta \boldsymbol{x}_{\mathrm{D},k}^{\mathrm{B}}(t) \end{bmatrix}. \tag{22}$$

Equation (22) reveals that increments of independent variables and dual variables can be represented parametrically as a linear function of boundary variables' increments. Next, substituting (22) into (21), and then simplifying it, the following quadratic function with respect to $\Delta \boldsymbol{x}_{\mathrm{D},k}^{\mathrm{B}}(t)$ can be obtained:

$$\boldsymbol{\Phi}_{\mathrm{D},k}(\Delta \boldsymbol{x}_{\mathrm{D},k}^{\mathrm{B}}(t)) = \frac{1}{2} \left( \Delta \boldsymbol{x}_{\mathrm{D},k}^{\mathrm{B}}(t) \right)^{\mathrm{T}} \boldsymbol{J}_{k,2} \Delta \boldsymbol{x}_{\mathrm{D},k}^{\mathrm{B}}(t) + \left( \boldsymbol{J}_{k,1} \right)^{\mathrm{T}} \Delta \boldsymbol{x}_{\mathrm{D},k}^{\mathrm{B}}(t) + \boldsymbol{J}_{k,0}, \tag{23}$$

where $J_{k,2}$ and $J_{k,1}$ denote the coefficient matrix and vector, respectively, corresponding to the quadratic and linear terms of the quadratic functions; and $J_{k,0}$ is constant. Thus, the quadratic function (23) can be used to describe the effect of the DS $k$ on the TS.

*3.2.2 Computation of boundary variables' increments*

Similarly, using the PDIPM, a sparse linear system can be constructed for the TS and convert it to the following quadratic programming problem $\Phi_T(\Delta x_T(t), \Delta y_T(t))$:

$$\min_{\Delta x_T(t), \Delta y_T(t)} \frac{1}{2}(\Delta x_T(t))^T H_T \Delta x_T(t) - (R_{T,x})^T \Delta x_T(t) + \Delta y_T(t)^T (\nabla_x g_T(x,t) - R_{T,y}) . \tag{24}$$

When receiving all quadratic functions with respect to boundary variables' increments from each DS, the variables' increments of the TS can be obtained by solving the following accumulated programming problem:

$$\{\Delta x_T^*(t), \Delta y_T^*(t)\} = \arg\min_{\Delta x_T(t), \Delta y_T(t)} \left\{ \Phi_T(\Delta x_T(t), \Delta y_T(t)) + \sum_{k \in \Omega_D} \Phi_{D,k}(\Delta x_{D,k}^B(t)) \right\} . \tag{25}$$

Next, as displayed in Fig. 4, after receiving $\Delta x_{D,k}^{B*}(t)$ from the TS, the DS $k$ can compute the increments of its independent variables and dual variables $\{\Delta x_{D,k}^{I*}(t), \Delta y_{D,k}^*(t)\}(k \in \Omega_D)$ by (22). The computed primal-dual variables' increments also satisfy the complete sparse linear system of the centralized model (2a)–(2e). The proof is in the appendix.

*3.2.3 Communication topology*

According to the variable partition structure in Fig. 3, all boundary variables exist in the TS model. Thus, a radial communication topology diagram can be designed to achieve a decentralized implementation for solving the sparse linear systems of the TS and DSs. Fig. 5 displays the interaction between the TS and DSs. When the TS receives the quadratic function from DSs, the boundary variables' increments can be obtained by (25) and transmitted to the DSs. Next, for each DS, after receiving $\Delta x_{D,k}^{B*}(t)$ from the TS, (22) can be used to compute its independent variables' increments. Fig. 4 indicates that only one exchange is conducted between the TS and DSs. Thus, such an implementation can considerably decrease the communication burden.

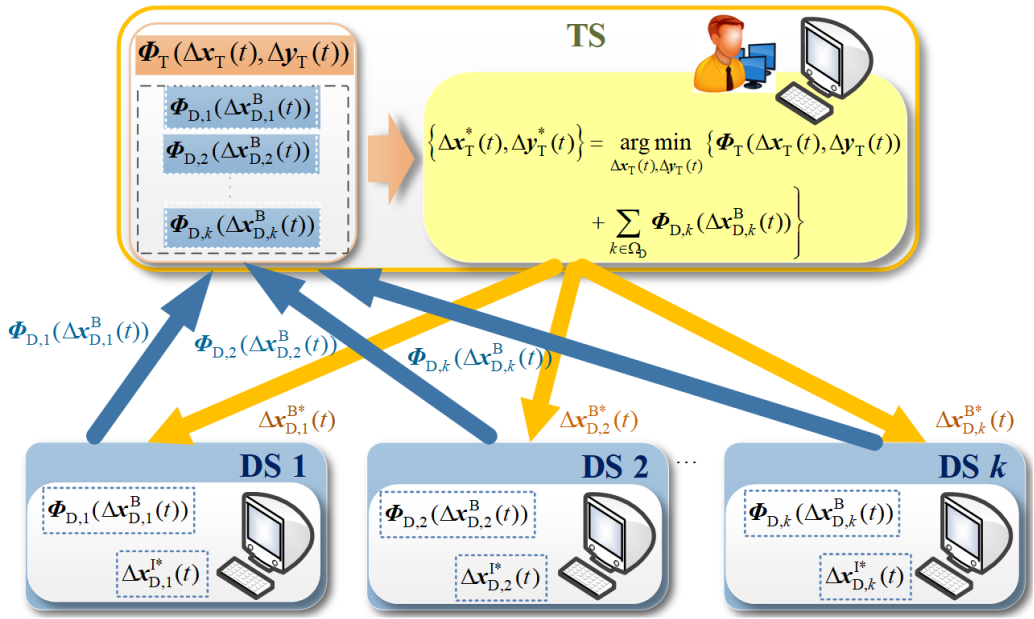

**Fig. 5. Radial communication graph and message passing between the TS and DSs.**

*3.2.4 Variable updating*

In the PDIPM, the scalar primal-dual step sizes $\{\alpha_p^{\mathrm{T}}, \alpha_d^{\mathrm{T}}\}, \{\alpha_p^k \alpha_d^k\}(k \in \Omega_{\mathrm{D}})$ are selected to preserve the feasible condition on slack variables $u$ and $l$, and Lagrange multipliers $w$ and $z$. Therefore, the variables of the dynamic system (12) are updated using the following variables:

$$x_{\mathrm{T}}(t+\tau) = x_{\mathrm{T}}(t) + \tau\alpha_p^{\mathrm{T}}\Delta x_{\mathrm{T}}(t), \tag{26a}$$

$$y_{\mathrm{T}}(t+\tau) = y_{\mathrm{T}}(t) + \tau\alpha_d^{\mathrm{T}}\Delta y_{\mathrm{T}}(t), \tag{26b}$$

$$l_{\mathrm{T}}(t+\tau) = l_{\mathrm{T}}(t) + \tau\alpha_p^{\mathrm{T}}\Delta l_{\mathrm{T}}(t), \tag{26c}$$

$$z_{\mathrm{T}}(t+\tau) = z_{\mathrm{T}}(t) + \tau\alpha_d^{\mathrm{T}}\Delta z_{\mathrm{T}}(t), \tag{26d}$$

$$u_{\mathrm{T}}(t+\tau) = u_{\mathrm{T}}(t) + \tau\alpha_p^{\mathrm{T}}\Delta u_{\mathrm{T}}(t), \tag{26e}$$

$$w_{\mathrm{T}}(t+\tau) = w_{\mathrm{T}}(t) + \tau\alpha_d^{\mathrm{T}}\Delta w_{\mathrm{T}}(t), \tag{26f}$$

and

$$x_{\mathrm{D},k}(t+\tau) = x_{\mathrm{D},k}(t) + \tau\alpha_p^k\Delta x_{\mathrm{D},k}(t), \ \forall k \in \Omega_{\mathrm{D}}, \tag{27a}$$

$$y_{\mathrm{D},k}(t+\tau) = y_{\mathrm{D},k}(t) + \tau\alpha_d^k\Delta y_{\mathrm{D},k}(t), \ \forall k \in \Omega_{\mathrm{D}}, \tag{27b}$$

$$l_{\mathrm{D},k}(t+\tau) = l_{\mathrm{D},k}(t) + \tau\alpha_p^k\Delta l_{\mathrm{D},k}(t), \ \forall k \in \Omega_{\mathrm{D}}, \tag{27c}$$

$$z_{\mathrm{D},k}(t+\tau) = z_{\mathrm{D},k}(t) + \tau\alpha_d^k\Delta z_{\mathrm{D},k}(t), \ \forall k \in \Omega_{\mathrm{D}}, \tag{27d}$$

$$u_{\mathrm{D},k}(t+\tau) = u_{\mathrm{D},k}(t) + \tau\alpha_p^k\Delta u_{\mathrm{D},k}(t), \ \forall k \in \Omega_{\mathrm{D}}, \tag{27e}$$

$$w_{\mathrm{D},k}(t+\tau) = w_{\mathrm{D},k}(t) + \tau\alpha_d^k\Delta w_{\mathrm{D},k}(t), \ \forall k \in \Omega_{\mathrm{D}}. \tag{27f}$$

*3.3. Online tracking framework for the proposed algorithm*

To achieve online applications, the output power setting points of the RESs should be given quickly after each sampling. In this case, only one or a few iterations are affordable. Fig. 6 depicts the sequence diagram of the online implementation for the proposed algorithm. Suppose that the proposed algorithm starts at an initial time $t_0$ and achieves online tracking at time $t_k$; thus, the update results for all time $t \geq t_k$ can be applied directly into the command level. To speed up the tracking process, the proposed algorithm is allowed to iterate from initial time $t_0$ until convergence at time $t_c$; and during the period $t_0 \sim t_c$, the model input parameters are assumed to be time-invariant. For all time $t \geq t_c$, the proposed algorithm only executes one iteration during each short sampling period, including one time of algebraic computing, algebraic updating, information exchanging, and command sending.

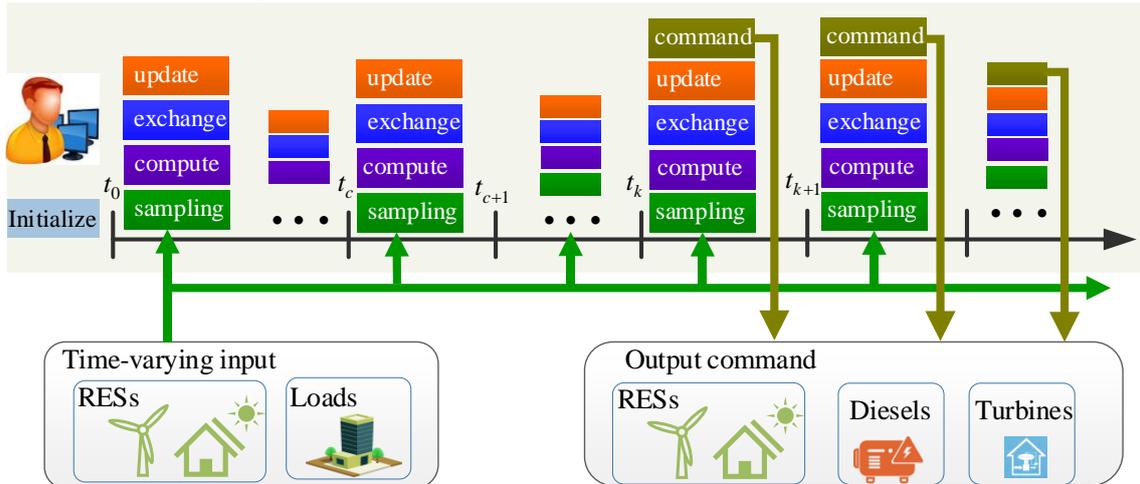

**Fig. 6. Sequence diagram of the proposed online distributed algorithm for online applications.**

## 4. Case studies

To analyze the performance of the proposed method, an integrated system was developed (Fig. 7). The IEEE 9-bus TS connects three DSs at buses 5, 7, and 9 of the TS. All the DSs are 33-bus systems, equipped with nine RESs, that is, four WT systems and five PV systems with the ratings of 200 kVA (Fig. 8).

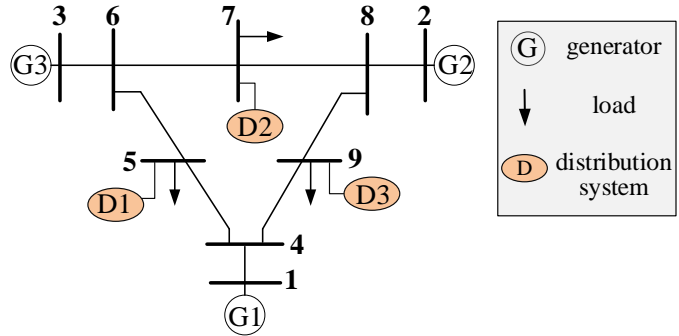

**Fig. 7. Topology of the IEEE 9-bus TS with three DSs.**

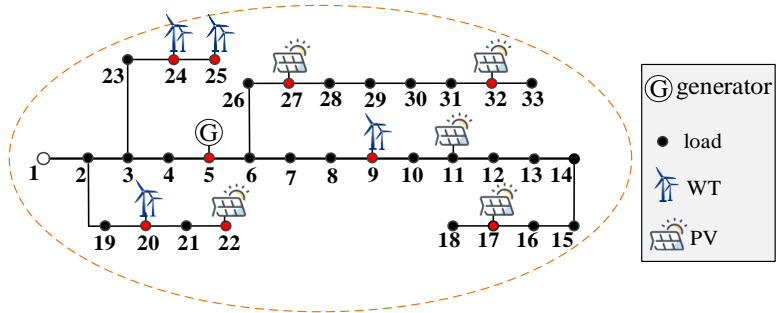

**Fig. 8. Topology of the 33-bus distribution system.**

To simulate a dynamic scenario, we used the real data of a particular day for PV/WT system generation and load profiles, available in [44]. Furthermore, we used Hermite interpolation to obtain continuously time-varying curves. Fig. 9 depicts the available active power of a PV system at bus 17 and a WT system at bus 25 as well as the total loads of the integrated system. We also diversify the load and PV/WT generation profile for each bus by small random additive noises. This study focused on online tracking for ultra-short-term (e.g., 20 ms) optimization, and because such short-term RES and load forecasts are typically accurate, their prediction errors are assumed to be negligible [45].

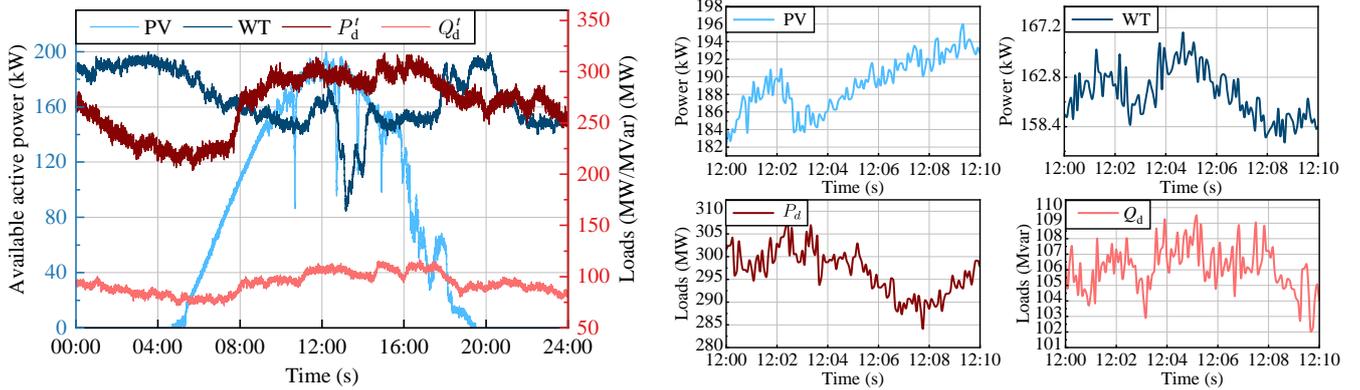

(a) One day  (b) Ten minutes at noon
**Fig. 9. Trajectories of PV and WT available active power and load.**

All the programs were developed in the MATLAB 2018b and run on a DELL Precision 3640 workstation with 64 GB RAM and 3.7 GHz Intel (R) Core (TM) i9-10900K CPU.

### 4.1. Coordination effectiveness analysis

To verify the correctness of the decentralized TD-ACOPF results, the centralized PDIPM [43] was used to obtain the centralized TD-ACOPF results for comparison with a sampling period of 20 ms. It is unpractical to obtain the centralized optimal trajectories in practical application because of privacy concerns and the limitation of computational capabilities. Next, we used the proposed decentralized approach to track the centralized optimal trajectories. Furthermore, the validity of coordination optimization was verified by using an independent optimization method by which the TS and DS optimize their own ACOPF based on a predefined boundary condition.

The trajectories of total objective values and voltage profiles of boundary buses are displayed in Fig. 10 and Fig. 11. The trajectories of the decentralized optimization are close to those of centralized optimization after approximately 1.5 s. This result implies that the proposed decentralized approach can achieve exact tracking for the optimum within a short time. Compared with the sampling period (i.e., 20 ms), the average computational time (i.e., 6 ms) and the maximum value (i.e., 9 ms) are almost trivial.

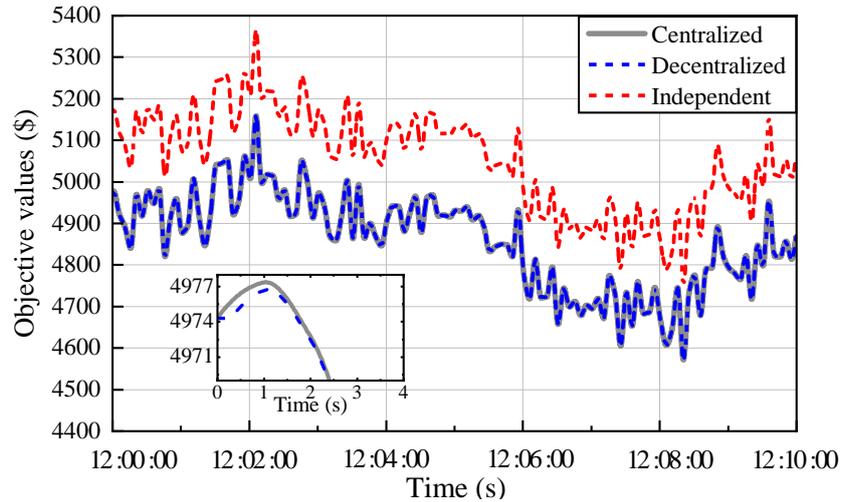

**Fig. 10. Total objective values obtained by centralized, decentralized, and independent optimization methods.**

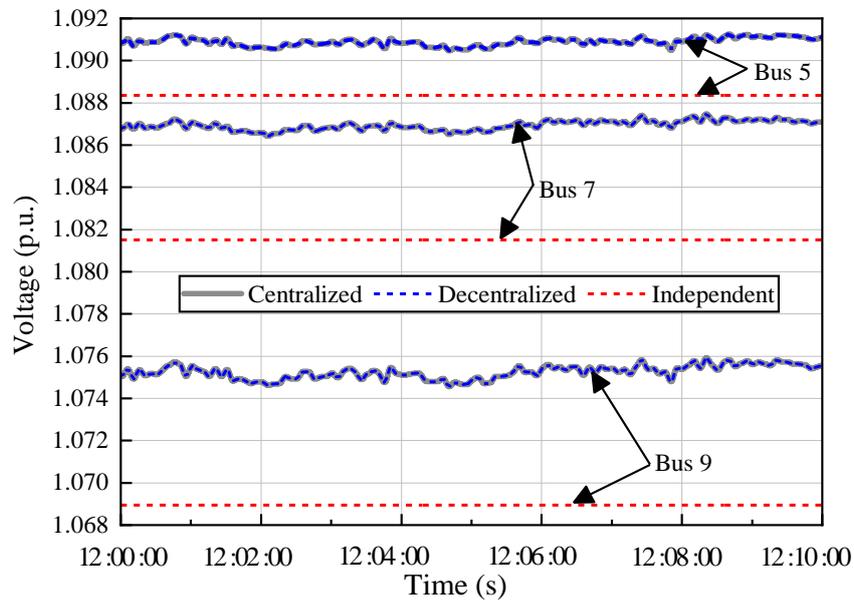

**Fig. 11. Voltage profiles of boundary buses (5, 7, and 9) obtained by centralized, decentralized, and independent optimization methods.**

As displayed in Fig. 10, the objective values of the decentralized optimization are always smaller than that of independent optimization, and the mean objective reduction ($195.11) accounts for 4.0% of the mean objective ($4857.38). Thus, the proposed decentralized optimization method can provide superior economic performance for the coupled transmission and distribution systems.

RES penetration can be increased three times. Fig. 12 presents a comparison of the voltage profiles of the distribution networks at $t$ = 12:02:00. Most of the bus voltages obtained by the independent optimization method reach the upper limit (1.1

p.u.), which likely leads to overvoltage in the presence of unexpected system disturbance. However, this risk can be mitigated by the proposed decentralized optimization because the narrow margin of voltage improved considerably, as displayed in Fig. 9. Thus, the proposed decentralized optimization method is highly suitable for handling the voltage rise issue of the distribution system.

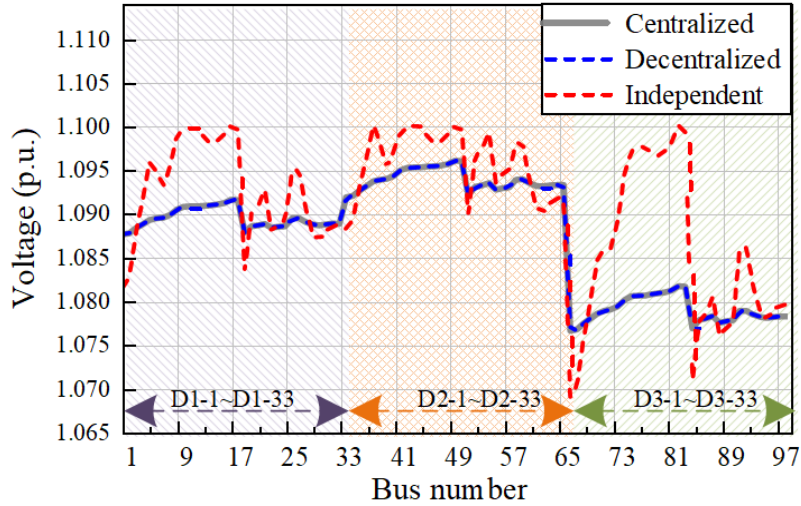

**Fig. 12. Voltage profiles obtained by centralized, decentralized, and independent optimization methods in case of high penetration of RESs.**

*4.2. Tracking performance analysis*

To analyze the effect of the prediction term, the proposed decentralized approach was modified by removing the prediction term in the dynamic system (15), and then the dynamic system without the prediction term was used to track centralized optimal trajectories. Fig. 13 displays the errors of the total objective values, which are defined as the relative deviations between the results obtained using the centralized method and those obtained using the decentralized methods with and without the prediction term. Using the prediction term in the proposed decentralized approach can decrease the tracking error.

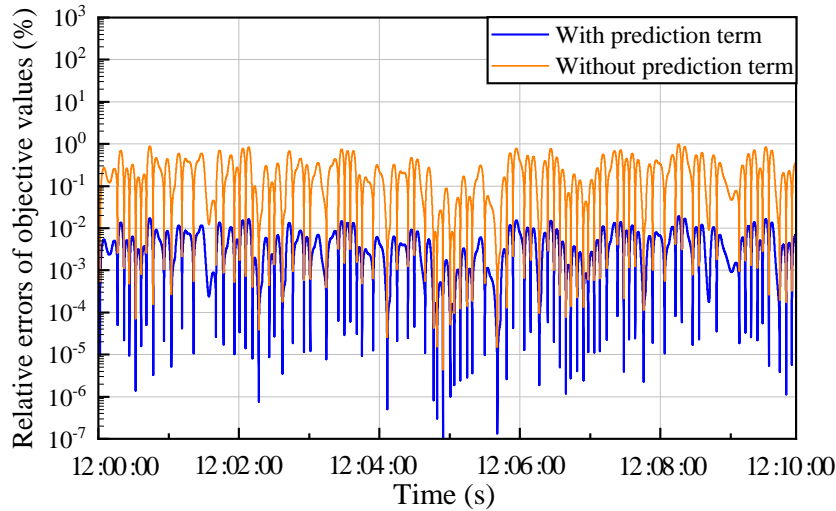

**Fig. 13. Relative errors of objective values obtained by the proposed decentralized approach with and without the prediction term.**

To verify the effect of the prediction term, we only depict the trajectories of objective values for 60 s from 12:00:00 to 12:01:00 in Fig. 14. Compared with the centralized results, the decentralized results without using the prediction term lag behind. Because the prediction term can be used to predict how the optimal solution changes over time, it can decrease the tracking error.

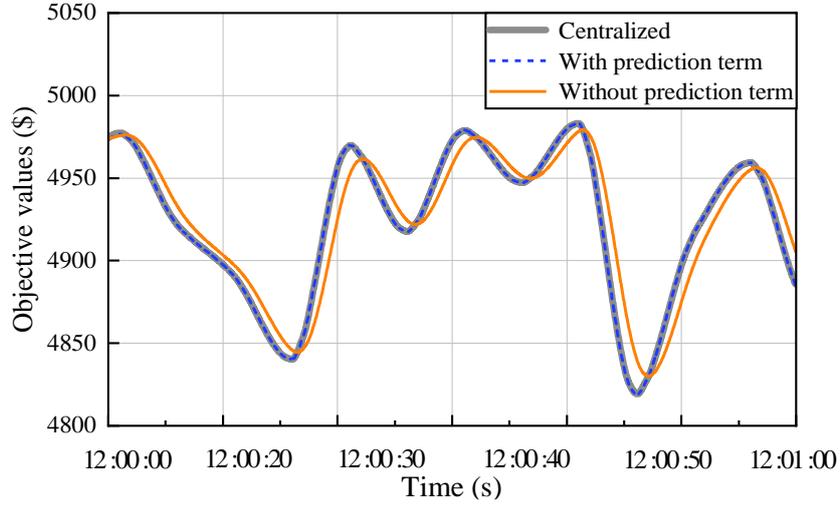

**Fig. 14. Objective values obtained by the centralized approach and proposed decentralized approach with and without the prediction term.**

To analyze the performance of the proposed decentralized approach for various sampling periods, four cases were simulated for $\Delta t \in \{0.01, 0.02, 0.05, 0.1, 0.5\}(s)$. As displayed in Fig. 15, the smaller the sampling period is, the smaller the tracking error is. Smaller periods denote smaller changes in adjacent periods, and the prediction is accurate. However, a too small sampling period (e.g., $\Delta t = 0.01s$) requires high communication and computation requirements. By contrast, a too large sampling period (e.g., $\Delta t = 0.5s$) results in poor tracking performance and even instability. Thus, we selected $\Delta t = 0.02s$ for compromise.

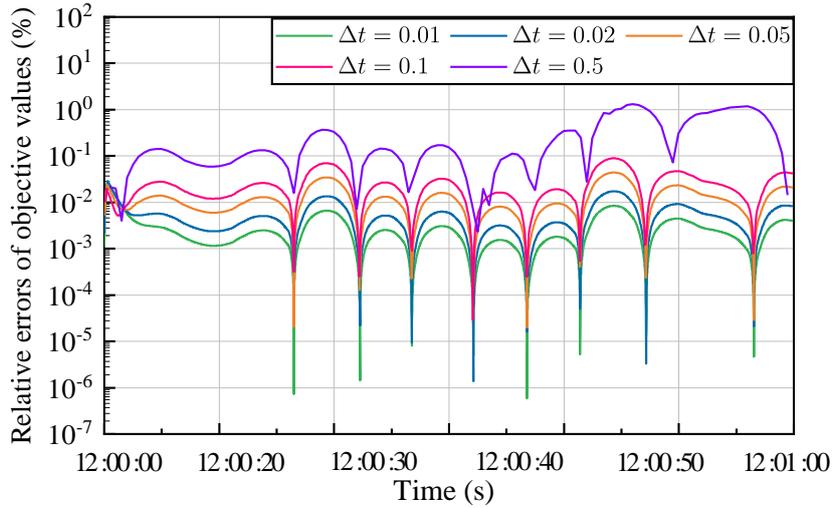

**Fig. 15. Relative errors of objective values obtained by the proposed decentralized approach with various sampling periods.**

*4.3. Online application analysis*

To analyze the online application capacity of the proposed decentralized method, the proposed decentralized approach with four conventional iterative-convergence algorithms, namely, HGD, OCD, ADMM, and APP, can be compared. The comparison terms include the computation time and iteration number in the range of 12:01:00 to 12:02:00 (total 60 s) in Table 1. The results indicate that the proposed algorithm only conducts one iteration and exhibits fast-tracking performance. However, the conventional iterative-convergence algorithm requires an excessive number of iterations, and the time consumption is unacceptable for online applications.

Table 1. Comparison of various algorithms for decentralized coordinated optimization.

|  | Iteration number | | Time (s) | |
| --- | --- | --- | --- | --- |
|  | Mean | Max | Mean | Max |
| Proposed | **1** | **1** | **0.006** | **0.009** |
| HGD | 35 | 48 | 1.98 | 2.38 |
| OCD | 59 | 74 | 3.01 | 5.80 |
| ADMM | 102 | 134 | 68.00 | 90.10 |
| APP | 97 | 128 | 61.20 | 88.00 |

(HGD [20]: heterogeneous decomposition; OCD [14]: optimality condition decomposition; ADMM [10]: alternating direction method of multipliers; APP [8]: the auxiliary problem principle)

## 5. Conclusion

An online decentralized tracking algorithm was proposed for the coordinated ACOPF in coupled transmission and distribution grids. The proposed method is designed from the perspective of control to construct a dynamic system to capture the information of the objectives and constraints that evolve over time. A decentralized implementation was developed based on only a few information exchanges with respect to boundary variables. The results of the case studies verified the coordination effectiveness of the proposed algorithm compared with the independent optimization method. The results also revealed that the performance of the proposed method is accurate, fast-tracking, and exhibits online application capacity that is greater than conventional iterative-convergence–based decentralized algorithms.

## 6. Appendix

Here, we prove that the computed primal-dual variables' increments by (28) also satisfy the centralized correction equation (28). On basis of the PDIPM, a sparse linear system (A-1) is formed and used to obtain the search direction $\Delta \hat{x}(t)$ and $\Delta \hat{y}(t)$, and then $\Delta \hat{w}(t), \Delta \hat{z}(t), \Delta \hat{u}(t)$, and $\Delta \hat{l}(t)$ (the superscript $\hat{}$ represents the union of variables from TS and DSs) can easily be calculated by (A-2).

$$\begin{bmatrix} \hat{H} & \nabla_{\hat{x}} \hat{g}(\hat{x},t) \\ \nabla_{\hat{x}}^T \hat{g}(\hat{x},t) & 0 \end{bmatrix} \begin{bmatrix} \Delta \hat{x}(t) \\ \Delta \hat{y}(t) \end{bmatrix} = \begin{bmatrix} \hat{R}_{\hat{x}} \\ \hat{R}_{\hat{x}} \end{bmatrix}, \quad (A\text{-}1)$$

$$\begin{cases} \Delta \hat{u}(t) = -(L_{\hat{w}} + L_{\hat{w}t}) - \nabla_{\hat{x}}^T h(\hat{x},t)\Delta \hat{x}(t) \\ \Delta \hat{l}(t) = (L_{\hat{z}} + L_{\hat{z}t}) + \nabla_{\hat{x}}^T h(\hat{x},t)\Delta \hat{x}(t) \\ \Delta \hat{w}(t) = -I_{\hat{U}}^{-1}(L_{\hat{u}} + L_{\hat{u}t}) - I_{\hat{U}}^{-1} I_{\hat{W}} \Delta \hat{u}(t) \\ \Delta \hat{z}(t) = -I_{\hat{L}}^{-1}(L_{\hat{l}} + L_{\hat{l}t}) - I_{\hat{L}}^{-1} I_{\hat{Z}} \Delta \hat{l}(t) \end{cases}. \quad (A\text{-}2)$$

Solving the sparse linear equations (A-1) is equivalent to solve the following quadratic programming problem:

$$\min_{\Delta \hat{x}(t), \Delta \hat{y}(t)} \frac{1}{2}(\Delta \hat{x}(t))^T \hat{H} \Delta x(t) - (\hat{R}_{\hat{x}})^T \Delta \hat{x}(t) + \Delta \hat{y}(t)^T (\nabla_x \hat{g}(\hat{x},t) \Delta \hat{x}(t) - \hat{R}_{\hat{y}}). \quad (A\text{-}3)$$

According to the superposition principle, we have

$$\begin{cases} \hat{H} = \sum_{k \in \Omega_D} (\bar{E}_k)^T H_{D,k} \bar{E}_k + (\bar{E}_T)^T H_T \bar{E}_T \\ \hat{R}_{\hat{x}} = \sum_{k \in \Omega_D} (\bar{E}_k)^T R_{k,x} + (\bar{E}_T)^T R_{T,x} \\ \hat{R}_{\hat{y}} = \sum_{k \in \Omega_D} (\bar{E}_k)^T R_{k,y} + (\bar{E}_T)^T R_{T,y} \end{cases}, \quad (A\text{-}4)$$

where $x_T(t) = \bar{E}_T \hat{x}(t), x_{D,k}(t) = \bar{E}_k \hat{x}(t), k \in \Omega_D$.

Substituting (A-4) into (A-3), the quadratic programming problem is rewritten as

$$\min_{\substack{\Delta x_T(t) \bigcup_{k \in \Omega_D} \Delta x_{D,k}(t), \\ \Delta y_T(t) \bigcup_{k \in \Omega_D} \Delta y_{D,k}(t)}} \left\{ \begin{aligned} & \frac{1}{2}(\Delta x_T(t))^T H_T \Delta x_T(t) - (R_{T,x})^T \Delta x_T(t) + \Delta y_T(t)^T (\nabla_{x_T} g_T(x_T,t)\Delta x_T(t) - R_{T,y}) \\ & + \sum_{k \in \Omega_D} \frac{1}{2}(\Delta x_{D,k}(t))^T H_{D,k} \Delta x_{D,k}(t) - (R_{k,x})^T \Delta x_{D,k}(t) + \Delta y_{D,k}(t)^T (\nabla_{x_{D,k}} g_{D,k}(x_T,t)\Delta x_{D,k}(t) - R_{k,y}) \end{aligned} \right\}. \quad (A\text{-}5)$$

According to the variables partition shown in Fig. 2, we rearrange $\Delta \boldsymbol{x}_\text{T}(t) \bigcup_{k \in \Omega_\text{D}} \Delta \boldsymbol{x}_{\text{D},k}(t)$ as $\left\{\Delta \boldsymbol{x}_\text{T}^\text{I}(t), \Delta \boldsymbol{x}^\text{B}, \bigcup_{k \in \Omega_\text{D}} \Delta \boldsymbol{x}_{\text{D},k}^\text{I}(t)\right\}$ and then the quadratic programming problem (A-5) can be equivalent to the following two-layer optimization form:

$$\min_{\Delta \boldsymbol{x}_\text{T}^\text{I}(t), \Delta \boldsymbol{x}^\text{B}(t)} \left\{ \begin{array}{l} \frac{1}{2}(\Delta \boldsymbol{x}_\text{T}(t))^\text{T} \boldsymbol{H}_\text{T} \Delta \boldsymbol{x}_\text{T}(t) - (\boldsymbol{R}_{\text{T},x})^\text{T} \Delta \boldsymbol{x}_\text{T}(t) + \Delta \boldsymbol{y}_\text{T}(t)^\text{T} (\nabla_{\boldsymbol{x}_\text{T}} \boldsymbol{g}_\text{T}(\boldsymbol{x}_\text{T},t) \Delta \boldsymbol{x}_\text{T}(t) - \boldsymbol{R}_{\text{T},y}) \\ + \min_{\bigcup_{k \in \Omega_\text{D}} \Delta \boldsymbol{x}_{\text{D},k}^\text{I}(t)} \sum_{k \in \Omega_\text{D}} \frac{1}{2} (\Delta \boldsymbol{x}_{\text{D},k}(t))^\text{T} \boldsymbol{H}_{\text{D},k} \Delta \boldsymbol{x}_{\text{D},k}(t) - (\boldsymbol{R}_{k,x})^\text{T} \Delta \boldsymbol{x}_{\text{D},k}(t) + \Delta \boldsymbol{y}_{\text{D},k}(t)^\text{T} (\nabla_{\boldsymbol{x}_{\text{D},k}} \boldsymbol{g}_{\text{D},k}(\boldsymbol{x}_\text{T},t) \Delta \boldsymbol{x}_{\text{D},k}(t) - \boldsymbol{R}_{k,y}) \end{array} \right\}. \quad \text{(A-6)}$$

Since the independent variable for each DS are independent, (A-6) can further be equivalent to the following form

$$\min_{\Delta \boldsymbol{x}_\text{T}^\text{I}(t), \Delta \boldsymbol{x}^\text{B}(t)} \left\{ \begin{array}{l} \frac{1}{2}(\Delta \boldsymbol{x}_\text{T}(t))^\text{T} \boldsymbol{H}_\text{T} \Delta \boldsymbol{x}_\text{T}(t) - (\boldsymbol{R}_{\text{T},x})^\text{T} \Delta \boldsymbol{x}_\text{T}(t) + \Delta \boldsymbol{y}_\text{T}(t)^\text{T} (\nabla_{\boldsymbol{x}_\text{T}} \boldsymbol{g}_\text{T}(\boldsymbol{x}_\text{T},t) \Delta \boldsymbol{x}_\text{T}(t) - \boldsymbol{R}_{\text{T},y}) \\ + \sum_{k \in \Omega_\text{D}} \min_{\Delta \boldsymbol{x}_{\text{D},k}^\text{I}(t)} \left( \frac{1}{2}(\Delta \boldsymbol{x}_{\text{D},k}(t))^\text{T} \boldsymbol{H}_{\text{D},k} \Delta \boldsymbol{x}_{\text{D},k}(t) - (\boldsymbol{R}_{k,x})^\text{T} \Delta \boldsymbol{x}_{\text{D},k}(t) + \Delta \boldsymbol{y}_{\text{D},k}(t)^\text{T} (\nabla_{\boldsymbol{x}_{\text{D},k}} \boldsymbol{g}_{\text{D},k}(\boldsymbol{x}_\text{T},t) \Delta \boldsymbol{x}_{\text{D},k}(t) - \boldsymbol{R}_{k,y}) \right) \end{array} \right\}, \quad \text{(A-7)}$$

From (A-7), it is clear that to solve the inner-layer optimization problem can form a quadratic function with respect to boundary variables' increments. Therefore, the outer-layer optimization problem is equivalent to (25).

## Acknowledgments

This work was supported in part by the National Natural Science Foundation of China under Grant No. 52107093, and in part by the Basic and Applied Basic Research Fund of Guangdong Province under Grant No. 2019A1515111037.

We thank LetPub (www.letpub.com) for its linguistic assistance during the preparation of this manuscript.